\documentclass[12pt]{amsart}
\usepackage{amsmath}
\usepackage{amssymb}
\usepackage{amscd}
\usepackage[mathscr]{eucal}
\usepackage{latexsym}

\newtheorem{thm}{Theorem}[section]
\newtheorem{lem}[thm]{Lemma}
\newtheorem{cor}[thm]{Corollary}      
\newtheorem{prop}[thm]{Proposition}

\theoremstyle{definition}

\theoremstyle{remark}

\numberwithin{equation}{section}

\def\M{{\mathcal M}}

\def\dR{\text{\rm dR}}
\def\Bri{\text{\rm Bri}}
\def\Sen{\text{\rm Sen}}
\def\cris{\text{\rm cris}}

\def\log{\text{\rm log}}
\def\dif{\text{\rm dif}}

\def\M{\text{\rm M}}
\def\st{\text{\rm st}}
\def\mod{\text{\rm mod}}

\def\Frac{\text{\rm Frac}}

\def\dim{\text{\rm dim}}
\def\-rig{\text{\rm -rig}}
\def\-log{\text{\rm -log}}
\def\-dif{\text{\rm -dif}}

\def\Gal{\text{\rm Gal}}

\def\Ker{\text{\rm Ker}\,}

\def\pf{\text{\rm pf}}

\def\lim{\text{\rm lim}}

\begin{document}

\title[Crystalline and semi-stable representations ]
{Crystalline and semi-stable representations in the imperfect residue field case}

\author[Kazuma Morita]{Kazuma Morita}
\address{Department of Mathematics, Hokkaido University, Sapporo 060-0810, Japan}
\email{morita@math.sci.hokudai.ac.jp}

\date{\today}

\maketitle

{\bf Abstract.}
Let $K$ be a $p$-adic local field with residue field $k$ such that $[k:k^{p}]=p^{e}<\infty$ and $V$ be a $p$-adic representation of $\Gal(\overline{K}/K)$. Then, by using the theory of $p$-adic differential modules, we show that $V$ is a potentially crystalline (resp. potentially semi-stable) representation of $\Gal(\overline{K}/K)$ if and only if $V$ is a potentially crystalline (resp. potentially semi-stable) representation of $\Gal(\overline{K^{\pf}}/K^{\pf})$ where $K^{\pf}/K$ is a certain $p$-adic local field whose residue field is the smallest perfect field $k^{\pf}$ containing $k$.   
As an application, we prove the $p$-adic monodromy theorem of Fontaine in the imperfect residue field case.

\section{Introduction}
Let $K$ be a complete discrete valuation field of characteristic $0$ with residue field $k$ of characteristic $p>0$ such that $[k:k^p]=p^e<\infty$.
Choose an algebraic closure $\overline{K}$ of $K$ and put $G_{K}=\Gal(\overline{K}/K)$.
By a $p$-adic representation of $G_{K}$, we mean a finite dimensional vector space $V$ over $\mathbb{Q}_{p}$ endowed with a continuous action of $G_{K}$.
As in the perfect residue field case, we can define the imperfect residue field versions of $B_{\cris}$ and $B_{\st}$ and, by using these rings, crystalline and semi-stable representations of $G_{K}$. 

Now, we shall state the main results of this article. Let us fix some notations.
Fix a lift $(b_{i})_{1\leq i\leq e}$ of a $p$-basis of $k$ in $\mathscr{O}_{K}$ (the ring of integers of $K$) and for each $m\geq 1$, fix a $p^{m}$-th root $b_{i}^{1/p^{m}}$ of $b_{i}$ in $\overline{K}$ satisfying $(b_{i}^{1/p^{m+1}})^{p}=b_{i}^{1/p^m}$.
Put $K^{(\pf)}=\cup_{m\geq 0}K(b_{i}^{1/p^{m}}, 1\leq i\leq e)$ and let $K^{\pf}$  be the $p$-adic completion of  $K^{(\pf)}.$
These fields depend on the choice of the sequences $(b_{i}^{1/p^{m}})_{m\in\mathbb{N}}$. Note that, if $V$ is a $p$-adic representation of $G_{K}$, it can be restricted to  a $p$-adic representation of $G_{K^{\pf}}=\Gal(\overline{K^{\pf}}/K^{\pf})$ where we choose an algebraic closure $\overline{K^{\pf}}$ of $K^{\pf}$ containing $\overline{K}$. 
Since $K^{\pf}$ is a complete discrete valuation field with perfect residue field, we can apply the classical theory (i.e. in the perfect residue field case)  to  $p$-adic representations of $G_{K^{\pf}}$.  
Our main results are the following.
\begin{thm} 
With  notation as above, we have the following equivalences.  
\begin{enumerate}
\item $V$ is a potentially crystalline representation of $G_{K}$ if and only if $V$ is a  potentially crystalline representation of     $G_{K^{\pf}}$,
\item $V$ is a  potentially semi-stable representation of $G_{K}$ if and only if $V$ is a  potentially semi-stable representation of $G_{K^{\pf}}$. 
\end{enumerate}
\end{thm} 
\begin{cor}
Keep the notation as in Theorem 1.1. Then, $V$ is a de Rham representation of $G_{K}$ if and only if $V$ is a potentially semi-stable representation of $G_{K}$.
\end{cor}
This paper is organized as follows.
In Section $2$, we shall review the definitions and basic known facts on crystalline and semi-stable representations, first in the perfect residue field case and then in the imperfect residue field case.
In Section $3$, we shall introduce some special elements which behave well under the action of $p$-adic differential operators. In Section $4$, by using these elements, we shall prove the main theorem. In Section $5$, as an application, we deduce the $p$-adic monodromy theorem of Fontaine  in the imperfect residue field case  (Corollary $1.2$) by using results of Berger [Be] and author [M]. 

{\bf Acknowledgments}
The author would like to thank Tetsushi Ito for useful discussions. He also would like to thank Marc-Hubert Nicole,  Shun Okubo and Seidai Yasuda for reading the manuscript carefully and giving useful comments.
He is grateful to his advisor Professor Kazuya Kato for his continuous advice and encouragements.
A part of this
work was done while the author was staying at Universit\'e Paris 11: he thanks this
institute for its hospitality. His staying at Universit\'e Paris 11 was partially
supported by JSPS Core-to-Core Program gNew Developments of Arithmetic Geometry,
Motives, Galois Theory, and Their Practical Applicationsh and he thanks
Professor Makoto Matsumoto for encouraging this visiting. This research was partially
supported by JSPS Research Fellowships for Young Scientists.  

\section{Review of crystalline and semi-stable representations}
\subsection{Crystalline and semi-stable representations in the perfect residue field case}
This subsection is a continuation of Subsection 2.1 of [M] and we keep the notation as in it. 
Let $\theta:\widetilde{\mathbb{A}}^{+}\rightarrow \mathscr{O}_{\mathbb{C}_{p}}$ be the natural homomorphism where $\mathscr{O}_{\mathbb{C}_{p}}$ denotes the ring of integers of $\mathbb{C}_{p}$. 
Define the ring $A_{\cris,K}$ to be the $p$-adic completion of the PD-envelope of $\Ker(\theta)$ compatible with the canonical PD-envelope over the ideal generated by $p$.
Put $B_{\cris,K}^{+}=A_{\cris,K}[1/p]$ and $B_{\cris,K}=B_{\cris,K}^{+}[1/t]$.
These rings are $K_{0}=\Frac(W(k))$-algebras endowed with an action of $G_{K}$ and an action of Frobenius $\varphi$ which commutes with the action of $G_{K}$. 
Furthermore, since we have the inclusion  $K\otimes_{K_{0}}B_{\cris,K}\hookrightarrow B_{\dR,K}$, the ring $K\otimes _{K_{0}}B_{\cris,K}$ is endowed with the filtration induced by that of $B_{\dR,K}$. 
Then, $(B_{\cris,K})^{G_{K}}$ is canonically isomorphic to $K_{0}$. Thus, for a $p$-adic representation $V$ of $G_{K}$, $D_{\cris,K}(V)=(B_{\cris,K}\otimes_{\mathbb{Q}_{p}}V)^{G_{K}}$ is naturally a $K_{0}$-vector space endowed with a Frobenius operator and a filtration after extending the scalars to $K$. We say that a $p$-adic representation $V$ of $G_{K}$ is a crystalline representation of $G_{K}$ if we have 
$$\dim_{\mathbb{Q}_p}V=\dim_{K_{0}}D_{\cris,K}(V)\quad \text{(we always have $\dim_{\mathbb{Q}_p}V\geq \dim_{K_{0}}D_{\cris,K}(V)$)}.$$
Furthermore, we say that a $p$-adic representation $V$ of $G_{K}$ is a potentially crystalline representation of $G_{K}$ if there exists a finite field extension $L/K$ in $\overline{K}$ such that $V$ is a crystalline representation of $G_{L}$.

Fix a prime element $\wp$ of $\mathscr{O}_{K}$ (the ring of integers of $K$) and an element $s=(s^{(n)})\in \widetilde{\mathbb{E}}^{+}$ such that $s^{(0)}=\wp$.
Then, the series $\log(s\wp^{-1})$ converges  to an element $u_{s}$ in $B_{\dR,K}^{+}$ and the subring  $B_{\cris,K}[u_{s}]$ of $B_{\dR,K}$ depends only on the choice of $\wp$.
We denote this ring by $B_{\st,K}$. Since we have the inclusion $K\otimes_{K_{0}}B_{\st,K}\hookrightarrow B_{\dR,K}$, the ring $K\otimes_{K_{0}}B_{\st,K}$ is  endowed with the action of $G_{K}$ and the filtration induced by that of $B_{\dR,K}$. 
The element $u_{s}$ is transcendental over $B_{\cris,K}$ and we extend the Frobenius $\varphi$ on  $B_{\cris,K}$ to  $B_{\st,K}$ by putting $\varphi(u_{s})=pu_{s}$. 
Furthermore, define the $B_{\cris,K}$-derivation $N:B_{\st,K}\rightarrow B_{\st,K}$ by $N(u_{s})=-1$.
It is easy to verify $N\varphi=p\varphi N$.
As in the case of $B_{\cris,K}$, we have $(B_{\st,K})^{G_{K}}=K_{0}$. Thus, for a $p$-adic representation $V$ of $G_{K}$,  $D_{\st,K}(V)=(B_{\st,K}\otimes_{\mathbb{Q}_{p}}V)^{G_{K}}$ is naturally a $K_{0}$-vector space endowed with a Frobenius operator and a filtration after extending the scalars to $K$. We say that a $p$-adic representation $V$ of $G_{K}$ is a semi-stable representation of $G_{K}$ if we have 
$$\dim_{\mathbb{Q}_p}V=\dim_{K_{0}}D_{\st,K}(V)\quad \text{(we always have $\dim_{\mathbb{Q}_p}V\geq \dim_{K_{0}}D_{\st,K}(V)$)}.$$
Furthermore, we say that a $p$-adic representation $V$ of $G_{K}$ is a potentially semi-stable representation of $G_{K}$ if there exists a finite field extension $L/K$ in $\overline{K}$ such that $V$ is a semi-stable representation of $G_{L}$.  
\subsection{Crystalline and semi-stable representations in the imperfect residue field case}
This subsection is a continuation of Subsection 2.2 of [M] and we keep the notation as in it.   
Let $k^{\pf}$ denote the perfect residue field of $K^{\pf}$ and put $K_{0}^{\pf}=\Frac( W(k^{\pf}))$.
Define $K_{0}=K\cap K_{0}^{\pf}$. Then, $K_{0}$ has residue field $k$ and the extension $K/K_{0}$ is finite.  Choose a Frobenius $\varphi$ on $K_{0}$ which is a lift of that on $k$. Put $\mathscr{O}_{K_{0}}=\mathscr{O}_{K}\cap W(k^{\pf})$.  
Let $\theta_{K_{0}}:\mathscr{O}_{K_{0}}\otimes_{\mathbb{Z}}\widetilde{\mathbb{A}}^{+}\rightarrow \mathscr{O}_{\mathbb{C}_{p}}$ denote the natural extension of $\theta:\widetilde{\mathbb{A}}^{+}\rightarrow \mathscr{O}_{\mathbb{C}_{p}}$ where    $\mathscr{O}_{\mathbb{C}_{p}}$ denotes the ring of integers of $\mathbb{C}_{p}$.
Define  $A_{\cris,K}$ to be the $p$-adic completion of the PD-envelope of $\Ker(\theta_{K_{0}})$ compatible with the canonical PD-envelope over the ideal generated by $p$.
Put $B_{\cris,K}^{+}=A_{\cris,K}[1/p]$ and $B_{\cris,K}=B_{\cris,K}^{+}[1/t]$.  
The ring $B_{\cris,K}$ is  the $K_{0}$-algebra endowed with an action of $G_{K}$ and an action of Frobenius $\varphi$ which commutes with the action of $G_{K}$.
Furthermore, since we have the inclusion  $K\otimes_{K_{0}}B_{\cris,K}\hookrightarrow B_{\dR,K}$ (see [Br, Proposition 2.47.]), the ring  $K\otimes_{K_{0}}B_{\cris,K}$ is endowed with the filtration induced by that of $B_{\dR,K}$.
For $1\leq i\leq e$, put $r_{i}=[\widetilde{b}_{i}]-b_{i}\in \mathscr{O}_{K_{0}}\otimes_{\mathbb{Z}} \widetilde{\mathbb{A}}^{+}$. Then, we have $r_{i}\in \Ker(\theta_{K_{0}})$ for $1\leq i\leq e$ and obtain an isomorphism
$$f:\text{$p$-adic completion of } \ A_{\cris,K^{\pf}}\langle r_{1},\ldots,r_{e}\rangle\rightarrow A_{\cris,K}$$
where $\langle *\rangle$ denotes PD-polynomial (see [Br, Proposition 2.39.]).
From this isomorphism, it follows that 
$$i:B_{\cris,K^{\pf}}\hookrightarrow B_{\cris,K} \quad\text{and}\quad p:B_{\cris,K}\twoheadrightarrow B_{\cris,K^{\pf}}:\mspace{5mu}r_{i}\mapsto 0$$
are $G_{K^{\pf}}$-equivariant homomorphisms and the composition 
$$p\circ i:B_{\cris,K^{\pf}}\hookrightarrow B_{\cris,K}\twoheadrightarrow B_{\cris,K^{\pf}}$$
is identity. By [Br, Proposition 2.50.], we have a canonical isomorphism $(B_{\cris,K})^{G_K}$ $=K_{0}$. Thus, for a $p$-adic representation $V$ of $G_{K}$, $D_{\cris,K}(V)=(B_{\cris,K}\otimes_{\mathbb{Q}_p}V)^{G_K}$ is naturally a $K_{0}$-vector space endowed with a Frobenius operator and a filtration after extending the scalars to $K$.
We say that a $p$-adic representation $V$ of $G_{K}$ is a crystalline representation of $G_{K}$ if we have
$$\dim_{\mathbb{Q}_p}V=\dim_{K_{0}}D_{\cris,K}(V).$$
Note that, for a $p$-adic representation $V$ of $G_{K}$, we always have $\dim_{\mathbb{Q}_p}V\geq \dim_{K_{0}}D_{\cris,K}(V)$ by [Br, Proposition 3.22.].  
Furthermore, we say that a $p$-adic representation $V$ of $G_{K}$ is a potentially crystalline representation of $G_{K}$ if there exists a finite field extension $L/K$ in $\overline{K}$ such that $V$ is a crystalline representation of $G_{L}$.

Fix a prime element $\wp$ of $\mathscr{O}_{K}$ and an element $s=(s^{(n)})\in \widetilde{\mathbb{E}}^{+}$ such that $s^{(0)}=\wp$.
Then, the series $\log(s\wp^{-1})$ converges  to an element $u_{s}$ in $B_{\dR,K}^{+}$ and the subring  $B_{\cris,K}[u_{s}]$ of $B_{\dR,K}$ depends only on the choice of $\wp$.
We denote this ring by $B_{\st,K}$. We can prove that  the element $u_{s}$ is transcendental over $B_{\cris,K}$ as in [F1, 4.3.]. Since we have the inclusion $K\otimes_{K_{0}}B_{\st,K}\hookrightarrow B_{\dR,K}$, the ring $K\otimes_{K_{0}}B_{\st,K}$ is  endowed with the action of $G_{K}$ and the filtration induced by that of $B_{\dR,K}$.  We extend the Frobenius $\varphi$ on $B_{\cris,K}$ to $B_{\st,K}$ by putting $\varphi(u_{s})=pu_{s}$. 
Furthermore, define the $B_{\cris,K}$-derivation $N:B_{\st,K}\rightarrow B_{\st,K}$ by $N(u_{s})=-1$.
It is easy to verify $N\varphi=p\varphi N$.
As in the case of $A_{\cris,K}$, we have an isomorphism 
$$f:(\text{$p$-adic completion of } \ A_{\cris,K^{\pf}} \langle r_{1},\ldots,r_{e}\rangle)[1/p,u_{s},1/t]\rightarrow B_{\st,K}$$
where $\langle *\rangle$ denotes PD-polynomial.
From this isomorphism, it follows that 
$$i:B_{\st,K^{\pf}}\hookrightarrow B_{\st,K} \quad\text{and}\quad p:B_{\st,K}\twoheadrightarrow B_{\st,K^{\pf}}:\mspace{5mu}r_{i}\mapsto 0$$
are $G_{K^{\pf}}$-equivariant homomorphisms and the composition 
$$p\circ i:B_{\st,K^{\pf}}\hookrightarrow B_{\st,K}\twoheadrightarrow B_{\st,K^{\pf}}$$
is identity.
By imitating the result [Br, Proposition 2.50.], we can show that we have a canonical isomorphism $(B_{\st,K})^{G_K}=K_{0}$. Thus, for a $p$-adic representation $V$ of $G_{K}$, $D_{\st,K}(V)=(B_{\st,K}\otimes_{\mathbb{Q}_p}V)^{G_K}$ is naturally a $K_{0}$-vector space endowed with a Frobenius operator and a filtration after extending the scalars to $K$. 
We say that a $p$-adic representation $V$ of $G_{K}$ is a semi-stable representation of $G_{K}$ if we have
$$\dim_{\mathbb{Q}_p}V=\dim_{K_{0}}D_{\st,K}(V).$$
For a $p$-adic representation $V$ of $G_{K}$, we can show that we always have $\dim_{\mathbb{Q}_p}V\geq \dim_{K_{0}}D_{\st,K}(V)$ in the same way as [Br, Proposition 3.22.].  
Furthermore, we say that a $p$-adic representation $V$ of $G_{K}$ is a potentially semi-stable representation of $G_{K}$ if there exists a finite field extension $L/K$ in $\overline{K}$ such that $V$ is a semi-stable representation of $G_{L}$.   
\section{Construction of special elements}
In this section, we shall introduce some special elements which behave well under the action of $p$-adic differential operators. Throughout this section, we keep the notation and assumptions of Section $3$ of [M]. For the definition of the $p$-adic differential modules $D_{\Sen}(V)$, $D_{\Bri}(V)$ and $D_{e\text{-}\dif}^{+}(V)$, refer to  $3.1.1$, $3.1.2$ and $3.1.3$ of [M]. Contrary to the definition of $D_{\dif}^{+}(V)$ of [M, 3.1.4], however,  define the module $D_{\dif}^{+}(V)$  as follows: for a $p$-adic representation $V$ of $G_{K^{\pf}}$, let $D_{\dif}^{+}(V)$ denote the union of $K_{\infty}^{\pf}[[t]]$-submodules of finite type of $(B_{\dR,K^{\pf}}^{+}\otimes_{\mathbb{Q}_{p}}V)^{H}$ stable under $\Gamma_{0}=\Gal(K_{\infty}^{\pf}/K^{\pf})$ (see [F2] for details).
\subsection{A special basis of $D_{e\text{-}\dif}^{+}(V)$}
We shall construct a special basis of $D_{e\text{-}\dif}^{+}(V)$ over $K_{\infty}^{(\pf)}[[t,t_{1},\ldots,t_{e}]]$ which bridges the gap between $D_{\dif}^{+}(V)$ and $D_{e\text{-}\dif}^{+}(V)$ and behaves well under the action of $\nabla^{(0)}$.  
Note that there is no $G_{K}$-equivariant injection $K\hookrightarrow B_{\dR,K^{\pf}}^{+}$:  
 we will sometimes write $L_{\dif}^{+}$ instead of the misleading $K_{\infty}^{\pf}[[t]]$.  
First, let us recall the following result.
\begin{prop} \verb+[+M, Proposition 4.8.\verb+]+ 
Let $V$ be a $p$-adic representation of $G_{K}$. If $V$ is a de Rham representation of $G_{K^{\pf}}$, there exists a $\nabla^{(0)}$-equivariant isomorphism of $K_{\infty}^{(\pf)}[[t,t_{1},\ldots,t_{e}]]$-modules
$$D_{e\text{-}\dif}^{+}(V)\simeq _{\nabla^{(0)}}\oplus_{j=1}^{d}K_{\infty}^{(\pf)}[[t,t_{1},\ldots,t_{e}]](n_{j})\quad (d=\dim_{\mathbb{Q}_{p}}V,\ n_{j}\in\mathbb{Z}).$$
\end{prop}
Next, let us define the $K_{\infty}^{\pf}[[t,t_{1},\ldots,t_{e}]]$-submodule $X$ of $(B_{\dR,K}^{+}\otimes_{\mathbb{Q}_{p}}V)^{H}$ by
$X=K_{\infty}^{\pf}[[t,t_{1},\ldots,t_{e}]]\otimes_{K_{\infty}^{(\pf)}[[t,t_{1},\ldots,t_{e}]]}D_{e\text{-}\dif}^{+}(V)$. If we put $D^{+,(r)}_{e\text{-}\dif}(V)=D^{+}_{e\text{-}\dif}(V)/(t,$ $t_{1},\ldots,t_{e})^{r}D^{+}_{e\text{-}\dif}(V)$, we have the inclusion $K_{\infty}^{\pf}\otimes_{K_{\infty}^{(\pf)}}D^{+,(r)}_{e\text{-}\dif}(V)\hookrightarrow L_{\dif}^{+}[[t_{1},\ldots,t_{e}]]/$ $(t,t_{1},\ldots,t_{e})^{r}\otimes _{L_{\dif}^{+}}D_{\dif}^{+}(V)$ by the theory of Sen. Since both sides have the same dimension over $K_{\infty}^{\pf}$, the inclusion above actually gives an isomorphism. By taking the projective limit with respect to $r$, we obtain a $\Gamma_{0}$-equivariant isomorphism 
$X\simeq L_{\dif}^{+}[[t_{1},\ldots,t_{e}]]\otimes_{L_{\dif}^{+}}D_{\dif}^{+}(V)$. 
\begin{prop} Let $V$ be a $p$-adic representation of $G_{K}$. If $V$ is a de Rham representation of $G_{K^{\pf}}$,  there exists a basis $\verb+{+f_{j}\verb+}+_{j=1}^{d}$ of $D_{\dif}^{+}(V)$ over $L_{\dif}^{+}$ such that 
\begin {enumerate}
\item $\verb+{+1\otimes f_{j}\verb+}+_{j=1}^{d}$ forms a basis of $D_{e\text{-}\dif}^{+}(V) \ (\subset X=L_{\dif}^{+}[[t_{1},\ldots,t_{e}]]\otimes_{L_{\dif}^{+}}D_{\dif}^{+}(V))$ over $K_{\infty}^{(\pf)}[[t,t_{1},\ldots,t_{e}]]$, 
\item the action of $\nabla^{(0)}$ on $\verb+{+1\otimes f_{j}\verb+}+_{j=1}^{d}$ is given by $\nabla^{(0)}(1\otimes f_{j})=n_{j}(1\otimes f_{j})$ 
 where the integers $n_{j}$ are those of Proposition $3.1.$
\end{enumerate}
\end{prop}
\begin{proof}
Let $\verb+{+G_{j}\verb+}+_{j=1}^{d}$ denote a basis of $D_{\dif}^{+}(V)$ over $K_{\infty}^{\pf}[[t]]$.
Since $D_{e\text{-}\dif}^{+}(V)$ is a submodule of $X=L_{\dif}^{+}[[t_{1},\ldots,t_{e}]]\otimes_{L_{\dif}^{+}}D_{\dif}^{+}(V)$, any element of $D_{e\text{-}\dif}^{+}(V)$ can be written as linear combinations of $\verb+{+1\otimes G_{j}\verb+}+_{j=1}^{d}$ over $L_{\dif}^{+}[[t_{1},\ldots,t_{e}]]$.
On the other hand, fix a basis $\verb+{+F_{j}\verb+}+_{j=1}^{d}$ of $D_{e\text{-}\dif}^{+}(V)$ over $K_{\infty}^{(\pf)}[[t,t_{1},\ldots,t_{e}]]$ that gives the isomorphism of Proposition $3.1$, that is, 
$\nabla^{(0)}(F_{j})=n_{j}F_{j}$ with  $n_{j}\in\mathbb{Z}$. 
Then,  we can write  
\begin{align}1\otimes F_{j}=\sum_{(m_{1},\ldots,m_{e})\in\mathbb{N}^{e}}t^{m_{1}}_{1}\cdots t_{e}^{m_{e}}\otimes (\sum_{k=1}^{d}a_{jk}^{(m_{1},\ldots,m_{e})}G_{k})\end{align}
where the $a_{jk}^{(m_{1},\ldots,m_{e})}$ are elements of $L_{\dif}^{+}$.
Put $f_{j}=\sum_{k=1}^{d}a_{jk}^{(0,\ldots,0)} G_{k}\in D_{\dif}^{+}(V).$  
Then, it follows that we have $\nabla^{(0)}(f_{j})=n_{j}f_{j}$. On the other hand, we have $\verb+{+\overline{f_{j}}=\overline{F_{j}}\verb+}+_{j=1}^{d}$ in $D_{\Sen}(V)$ where $-$ denotes the reduction modulo $(t,t_{1},\ldots,t_{e})X$. Since $\verb+{+\overline{F_{j}}\verb+}+_{j=1}^{d}$ forms a basis of $D_{\Sen}(V)$ over $K_{\infty}^{\pf}$, 
the lift $\verb+{+1\otimes f_{j}\verb+}+_{j=1}^{d}$ of $\verb+{+\overline{f_{j}}=\overline{F_{j}}\verb+}+_{j=1}^{d}$ in $X$ forms a basis of $X$ over $K_{\infty}^{\pf}[[t,t_{1},\ldots,t_{e}]]$. Furthermore, since $\verb+{+f_{j}\verb+}+_{j=1}^{d}$ are elements of $D_{\dif}^{+}(V)$, it follows that $\verb+{+f_{j}\verb+}+_{j=1}^{d}$ also forms a basis of $D_{\dif}^{+}(V)$ over $K_{\infty}^{\pf}[[t]]$. Thus, it remains to show that $\verb+{+1\otimes f_{j}\verb+}+_{j=1}^{d}$ forms a basis of $D_{e\text{-}\dif}^{+}(V)$ over $K_{\infty}^{(\pf)}[[t,t_{1},\ldots,t_{e}]]$.  Put $X_{r}=X/(t,t_{1},\ldots,t_{e})^{r}X$. Let $Y_{r}$ denote the $K_{\infty}^{(\pf)}[[t,t_{1},\ldots,t_{e}]]$-submodule of $X_{r}$ generated by the finite set $\verb+{+\sum_{k=1}^{d}a_{jk}^{(m_{1},\ldots,m_{e})} G_{k}\verb+}+_{j,m_{1}+\cdots  +m_{e}<r}\subset  (B_{\dR,K^{\pf}}^{+}\otimes_{\mathbb{Q}_{p}}V)^{H}$.  Then, it follows that this finitely generated $K_{\infty}^{(\pf)}[[t,t_{1},\ldots, t_{e}]]$-module $Y_{r}$ is stable under the action of $\Gamma_{K}$ by ($3.1$) and thus is contained in $D_{e\text{-}\dif}^{+,(r)}(V)$ by definition. On the other hand, $Y_{r}$ contains the elements $\verb+{+1\otimes f_{j}\verb+}+_{j=1}^{d}$ which are linearly independent over $K_{\infty}^{\pf}[[t,$ $t_{1},\ldots,t_{e}]]/(t,t_{1},\ldots,t_{e})^{r}$. Thus, both of $Y_{r}$ and $D_{e\text{-}\dif}^{+,(r)}(V)$ have the same dimension over $K_{\infty}^{(\pf)}$ and we get the equality $Y_{r}=D_{e\text{-}\dif}^{+,(r)}(V)$.
Therefore, by taking the projective limit with respect to $r$, we conclude that $\verb+{+1\otimes f_{j}\verb+}+_{j=1}^{d} \ (\subset \varprojlim_{r}Y_{r})$ forms a basis of $D_{e\text{-}\dif}^{+}(V)$ over $K_{\infty}^{(\pf)}[[t,t_{1},\ldots,t_{e}]]$.
\end{proof}
\begin{lem}By restricting $\nabla^{(i)}:D_{e\text{-}\dif}^{+}(V)\rightarrow D_{e\text{-}\dif}^{+}(V)$ $(0\leq i\leq e)$, we obtain  $\nabla^{(i)}:D_{e\text{-}\dif}^{+}(V)\cap (B_{\dR,K^{\pf}}^{+}\otimes_{\mathbb{Q}_{p}}V)^{H}\rightarrow D_{e\text{-}\dif}^{+}(V)\cap (B_{\dR,K^{\pf}}^{+}\otimes_{\mathbb{Q}_{p}}V)^{H}$ in  $(B_{\dR,K}^{+}\otimes_{\mathbb{Q}_{p}}V)^{H}$.
\end{lem}
\begin{proof}
For simplicity, put  $L_{\dR}^{+}=(B_{\dR,K^{\pf}}^{+})^{H}$, $ L_{\dR}^{+}(V)=(B_{\dR,K^{\pf}}^{+}\otimes_{\mathbb{Q}_{p}}V)^{H}$ and $Z=(B_{\dR,K}^{+}\otimes_{\mathbb{Q}_{p}}V)^{H}$.
Let $m_{\dR}$ denote the maximal ideal $(t,t_{1},\ldots,t_{e})$ of $(B_{\dR,K}^{+})^{H}$.
Then, we have 
\begin{align*}Z=\varprojlim \text{$_{r}$}Z/&m_{\dR}^{r}Z \quad \supset \quad L_{\dR}^{+}(V)=\varprojlim\text{$_{r}$}L_{\dR}^{+}(V)/(m_{\dR}^{r}Z\cap L_{\dR}^{+}(V))\\
\cup&\\
D_{e\text{-}\dif}^{+}(V)=&\varprojlim \text{$_{r}$}D_{e\text{-}\dif}^{+}(V)/(m_{\dR}^{r}Z\cap D_{e\text{-}\dif}^{+}(V)).
\end{align*}
Define $W$ as  the $L_{\dR}^{+}\cap K_{\infty}^{(\pf)}[[t,t_{1},\ldots, t_{e}]]$-submodule of $Z$ generated by  $L_{\dR}^{+}(V)\cap D_{e\text{-}\dif}^{+}(V)$. If we put $\hat{W}=\varprojlim_{r}W_{r}$ where $W_{r}$ denotes $W/(m_{\dR}^{r}Z\cap W)$, we have $L_{\dR}^{+}(V)\supset \hat{W}$ and $D_{e\text{-}\dif}^{+}(V)\supset \hat{W}$. Thus, we obtain $\hat{W}=W$ by definition. Therefore, it suffices to show that $W_{r}$ is stable under the actions of $\verb+{+\nabla^{(i)}\verb+}+_{i=0}^{e}$.  
Fix a basis $\verb+{+g_{j}\verb+}+_{j=1}^{h}$ of $D_{e\text{-}\dif}^{+}(V)/(m_{\dR}^{r}Z\cap D_{e\text{-}\dif}^{+}(V))$ over $K_{\infty}^{(\pf)}$. 
Then, there exists a finite field extension $L/K$ in $K_{\infty}^{(\pf)}$ such that  $\oplus_{j=1}^{h}L\cdot g_{j}$ is stable by the action of $\Gamma_{K}=G_{K}/H=\Gal(K_{\infty}^{(\pf)}/K)$. 
Thus, there exists an open subgroup $\Gamma_{i}'$ of $\Gamma_{i}$ such that, for all $\gamma\in\Gamma_{0}'$ (resp. $\beta_{i}\in\Gamma_{i}'$), the action of $\gamma$ (resp. $\beta_{i}$) on $\oplus_{j=1}^{h}L\cdot g_{j}$ is $L$-linear. Then, the series 
$$\log(\gamma)=-\sum_{k= 1}^{\infty}\frac{(\gamma-1)^{k}}{k} \quad (\text{resp.}\  \log(\beta_{i})=-\sum_{k= 1}^{\infty}\frac{(\beta_{i}-1)^{k}}{k})$$ converges to an endomorphism  of $\oplus_{j=1}^{h}L\cdot g_{j}.$   
These actions on $\oplus_{j=1}^{h}L\cdot g_{j}$ can be extended to those on  $D_{e\text{-}\dif}^{+}(V)/(m_{\dR}^{r}Z\cap D_{e\text{-}\dif}^{+}(V))$ by $K_{\infty}^{(\pf)}$-linearity. Since $W_{r}$ is contained in $D_{e\text{-}\dif}^{+}(V)/(m_{\dR}^{r}Z\cap D_{e\text{-}\dif}^{+}(V))$ and  stable under the action of $\Gamma_{K}$, it follows that $W_{r}$ is equipped with actions of $\nabla^{(0)}=\frac{\log(\gamma)}{\log(\chi(\gamma))}$ and $\nabla^{(i)}=\frac{\log(\beta_{i})}{c_{i}(\beta_{i})}$.
\end{proof}
\subsection{ $\widetilde{D}_{\cris,K^{\pf}}(V)$ and $\widetilde{D}_{\st,K^{\pf}}(V)$}
In this subsection, for simplicity, we shall denote $\widetilde{B}_{\cris,K^{\pf}}=(B_{\cris,K^{\pf}})^{H}$  and 
$\widetilde{D}_{\cris,K^{\pf}}(V)=(B_{\cris,K^{\pf}}\otimes_{\mathbb{Q}_{p}}V)^{H}$ (resp.  $\widetilde{B}_{\st,K^{\pf}}=(B_{\st,K^{\pf}})^{H}$ and  $\widetilde{D}_{\st,K^{\pf}}(V)=(B_{\st,K^{\pf}}\otimes_{\mathbb{Q}_{p}}V)^{H}).$
\begin{prop}(cf. Proposition $3.2$.)
Let $V$ be a $p$-adic representation of $G_{K}$. If $V$ is a crystalline (resp. semi-stable) representation of $G_{K^{\pf}}$, there exists a basis $\verb+{+g_{j}\verb+}+_{j=1}^{d}$ of $\widetilde{D}_{\cris,K^{\pf}}(V)$ over $\widetilde{B}_{\cris,K^{\pf}}$ (resp.  $\widetilde{D}_{\st,K^{\pf}}(V)$ over $\widetilde{B}_{\st,K^{\pf}}$) such that 
\begin{enumerate}
\item $\verb+{+g_{j}\verb+}+_{j=1}^{d}$ forms a basis of $D_{e\text{-}\dif}^{+}(V)[1/t]$ over $K_{\infty}^{(\pf)}[[t,t_{1},$ $\ldots,t_{e}]][1/t]$, 
\item  $\verb+{+g_{j}\verb+}+_{j=1}^{d}$ is contained in $\Ker(\nabla^{(0)})\ (\subset D_{e\text{-}\dif}^{+}(V)[1/t])$.   
\end{enumerate}
\end{prop}
\begin{proof}
Since the semi-stable representation case is similar, we shall consider only
the crystalline representation case. 
Since $V$ is also a de Rham representation of $G_{K^{\pf}}$, 
by Proposition $3.2$, there exists a basis $\verb+{+f_{j}\verb+}+_{j=1}^{d}$ of $D_{\dif}^{+}(V)$ over $K_{\infty}^{\pf}[[t]]$  such that (1) $\verb+{+1\otimes f_{j}\verb+}+_{j=1}^{d}$  forms a basis of $D_{e\text{-}\dif}^{+}(V)$ over $K_{\infty}^{(\pf)}[[t,t_{1},\ldots,t_{e}]]$ and (2) $\nabla^{(0)}(1\otimes f_{j})=n_{j}(1\otimes f_{j})$ with $n_{j}\in\mathbb{Z}$.
Then, since the action of $\nabla^{(0)}$ on $\verb+{+f_{j}t^{-n_{j}}\verb+}+_{j=1}^{d}$ is trivial and $\verb+{+f_{j}t^{-n_{j}}\verb+}+_{j=1}^{d}$ is contained in $D_{\dif}^{+}(V)[1/t]\subset (B_{\dR,K^{\pf}}\otimes V)^{H}$, there exists a finite field extension $L^{\pf}/K^{\pf}$ in $K_{\infty}^{\pf}$ such that $\verb+{+f_{j}t^{-n_{j}}\verb+}+_{j=1}^{d}$ forms a basis of $D_{\dR,L^{\pf}}(V)$ over $L^{\pf}$. If  $K=K_{0}(\alpha)$ and $L^{\pf}=K^{\pf}(\beta)$ for some $\beta=\zeta_{p^{n}}\in K_{\infty}^{\pf}$, there exists an element $a\in K_{\infty}^{(\pf)}$ such that $K_{0}(\alpha,\beta)=K_{0}(a)$. Then, we have $L^{\pf}=K_{0}^{\pf}(\alpha,\beta)=K_{0}^{\pf}(a)=L^{\pf}_{0}(a)$. Since $V$ is also a crystalline representation of $G_{L^{\pf}}$, we have  $D_{\dR,L^{\pf}}(V)=L^{\pf}_{0}(a)\otimes _{L_{0}^{\pf}}D_{\cris, L^{\pf}}(V)$. Thus, we can write 
\begin{align}f_{j}t^{-n_{j}}=\sum_{k=0}^{\delta-1}a^{k}\otimes g_{jk}\quad(g_{jk}\in D_{\cris, L^{\pf}}(V), \ \delta=[L_{0}^{\pf}(a):L_{0}^{\pf}]).\end{align}
We can extract a basis of $D_{\cris, L^{\pf}}(V)$ over $L_{0}^{\pf}$ from the family $\verb+{+g_{jk}\verb+}+_{j,k}$:  denote it by $\verb+{+g_{j}\verb+}+_{j=1}^{d}$. Since we have $B_{\cris,K^{\pf}}\otimes_{L_{0}^{\pf}}D_{\cris, L^{\pf}}(V)\simeq  B_{\cris,K^{\pf}}\otimes_{\mathbb{Q}_{p}}V$, by taking the invariant part by $H$, it follows that $\verb+{+g_{j}\verb+}+_{j=1}^{d}$ forms a basis of $\widetilde{D}_{\cris,K^{\pf}}(V)$ over $\widetilde{B}_{\cris,K^{\pf}}$. 
Furthermore, by ($3.2$), the action of $\nabla^{(0)}$ on $\verb+{+g_{j}\verb+}+_{j=1}^{d}$ is trivial. Thus, it remains to show that 
$\verb+{+ g_{j}\verb+}+_{j=1}^{d}$ forms a basis of $D_{e\text{-}\dif}^{+}(V)[1/t]$ over $K_{\infty}^{(\pf)}[[t,t_{1},\ldots,t_{e}]][1/t]$.  
First, let $Z_{r}$ denote the union of $K^{(\pf)}_{\infty}[[t,t_1,\ldots,t_e]]$-submodules 
of finite type of $(B_{\dR,K}^{+}\otimes_{\mathbb{Q}_{p}}V)^{H}/(t,t_{1},\ldots,t_{e})^{r}(B_{\dR,K}^{+}\otimes_{\mathbb{Q}_{p}}V)^{H}$ that are stable under the action of an open subgroup $\Gamma$ of $\Gamma_{K}$. Since we have the inclusion $D_{e\text{-}\dif}^{+,(r)}(V) \hookrightarrow Z_{r}$ by definition and both sides have the same dimension over $K^{(\pf)}_{\infty}$, we have  $D_{e\text{-}\dif}^{+,(r)}(V)=Z_{r}$. Thus, by taking the projective limit with respect to $r$, we obtain $D_{e\text{-}\dif}^{+}(V)=\varprojlim_{r}Z_{r}$. Choose integers $\verb+{+m_{jk}\verb+}+_{ 1\leq j \leq d, \  0\leq  k\leq  \delta-1}\subset \mathbb{Z}$ such that we have 
$$\verb+{+t^{m_{jk}} a^{k}\otimes g_{jk}\verb+}+_{ 1\leq j\leq d, \  0\leq k\leq \delta-1}\subset (B_{\dR,K}^{+}\otimes_{\mathbb{Q}_{p}}V)^{H}.$$ 
Let $Z$ denote the $K_{\infty}^{(\pf)}[[t,t_{1},\ldots,t_{e}]]$-submodule of $(B_{\dR,K}^{+}\otimes_{\mathbb{Q}_{p}}V)^{H}$ generated by the finite set $\verb+{+t^{m_{jk}} a^{k}\otimes g_{jk}\verb+}+_{1\leq j\leq d, \ 0\leq k\leq \delta-1}$. Take an open subgroup $\Gamma$ of $\Gamma_{K}$ such that the action of $\Gamma$ on the finite set $\verb+{+a^{k}\verb+}+_{k=0}^{\delta-1}$ is trivial. Then it follows from $(3.2)$ that this finitely generated $K_{\infty}^{(\pf)}[[t,t_{1},\ldots,t_{e}]]$-module $Z$ is stable under the action of $\Gamma$  and thus is contained in $D_{e\text{-}\dif}^{+}(V)$ by the preceding argument.  In particular, it follows that the elements  $\verb+{+g_{j}\verb+}+_{j=1}^{d}$ are contained in  $D_{e\text{-}\dif}^{+}(V)[1/t]$. 
Furthermore, since $\verb+{+g_{j}\verb+}+_{j=1}^{d}$ forms a basis of $B_{\dR,K}\otimes_{\mathbb{Q}_{p}}V$ over $B_{\dR,K}$, it is, in particular,  linearly independent over $B_{\dR,K}$ in $B_{\dR,K}\otimes_{\mathbb{Q}_{p}}V$.  
Take $m_{j}\in \mathbb{Z}$ such that we have 
$$\verb+{+g_{j}'=t^{m_{j}}g_{j}\verb+}+_{j=1}^{d}\subset  D_{e\text{-}\dif}^{+}(V).$$  
Let $\mathbb{L}(V)$ be the submodule of $B_{\dR,K}^{+}\otimes_{\mathbb{Q}_{p}}V$ generated by $\verb+{+g_{j}'\verb+}+_{j=1}^{d}$ over $B_{\dR,K}^{+}$ 
and let $\mathbb{D}(V)$ denote the union of $K_{\infty}^{(\pf)}[[t,t_{1},\ldots,t_{e}]]$-submodules of finite type of $\mathbb{L}(V)^{H}$ stable under 
$\Gamma_{K}$.  Since  $\verb+{+g_{j}'\verb+}+_{j=1}^{d} (\subset \mathbb{D}(V))$ forms a basis of $\mathbb{L}(V)$ over $B_{\dR,K}^{+}$, it follows that  $\verb+{+g_{j}'\verb+}+_{j=1}^{d}$  also forms a 
basis of $\mathbb{D}(V)$ over  $K_{\infty}^{(\pf)}[[t,t_{1},\ldots,t_{e}]]$ (see [A-B, Lemma 5.10]). For any element $x\in  D_{e\text{-}\dif}^{+}(V)[1/t]$, one can see that there exists an integer  $m\in \mathbb{Z}$ such that we have $t^{m}x\in  \mathbb{D}(V)$. Thus, $t^{m}x$ can be written as linear combinations of $\verb+{+g_{j}'\verb+}+_{j=1}^{d}$ over $K_{\infty}^{(\pf)}[[t,t_{1},\ldots,t_{e}]]$. It follows that   $\verb+{+g_{j}\verb+}+_{j=1}^{d}$ forms a basis of  $D_{e\text{-}\dif}^{+}(V)[1/t]$  over $K_{\infty}^{(\pf)}[[t,t_{1},\ldots,t_{e}]][1/t]$. 
\end{proof}
From now on, we shall keep the notation and assumptions of Proposition $3.4$.
The following result is proved in Proposition $3.5$ and Corollary $3.6$ of [M].
\begin{prop}
The action of $\verb+{+\nabla^{(i)}\verb+}+_{i=1}^{e}$ on the basis $\verb+{+g_{j}\verb+}+_{j=1}^{d}$ is given by
$(\nabla^{(1)})^{k_{1}}$ $\cdots(\nabla^{(e)})^{k_{e}}(g_{j})=t^{k_{1}+\cdots +k_{e}}\sum_{l=1}^{d}c_{j,\underline{k},l}g_{l}$ 
 where the $c_{j,\underline{k},l}$ $(\underline{k}=(k_{1},\ldots,k_{e}))$ are elements of $K_{\infty}^{(\pf)}[[t,t_{1},\ldots,t_{e}]]$ such that $\nabla^{(0)}(c_{j,\underline{k},l})=0$.
\end{prop}
\begin{prop}Let $V$ be a $p$-adic representation of $G_{K}$. If $V$ is  a crystalline (resp. semi-stable) representation of $G_{K^{\pf}}$, we have 
$$(\nabla^{(1)})^{k_{1}}\cdots(\nabla^{(e)})^{k_{e}}(g_{j})\in \widetilde{D}_{\cris,K^{\pf}}(V)\quad (\text{resp.} \ \in \widetilde{D}_{\st,K^{\pf}}(V))$$  
for all $(k_{i})_{1\leq i\leq e}\in\mathbb{N}^{e}$ and $1\leq j\leq d$.
\end{prop}
\begin{proof}Since  the semi-stable representation case is similar, we shall consider only the crystalline representation case.  
It is enough to prove that if $g\in D_{e\text{-}\dif}^{+}(V)[1/t]$ is such that $g\in \widetilde{D}_{\cris,K^{\pf}}(V)$ and  $\nabla^{(0)}(g)=0$, then $(\nabla^{(1)})^{k_{1}}\cdots(\nabla^{(e)})^{k_{e}}(g)\in \widetilde{D}_{\cris,K^{\pf}}(V)$  for all $(k_{i})_{1\leq i\leq e}\in\mathbb{N}^{e}$. Since the proof of the general case is exactly the same (only with heavier notations), we just show that $\nabla^{(i)}(g)\in \widetilde{D}_{\cris,K^{\pf}}(V)$ for $1\leq i\leq e$.
First, for $r\in\mathbb{N}_{>0}$ and $h\in D_{e\text{-}\dif}^{+}(V)$, there exists an open subgroup $\Gamma_{i}^{h,r}$ of $\Gamma_{i}$ such that we have $\beta_{i}(h)=\exp(c_{i}(\beta_{i})\nabla^{(i)})(h)$ mod $(t,t_{1},\ldots,t_{e})^{r}$ $D_{e\text{-}\dif}^{+}(V)$ for all $\beta_{i}\in \Gamma_{i}^{h,r}$ (see [A-B] and [F2]). Thus, if we take   $M\in\mathbb{N}$ such that $t^{M}g\in D_{e\text{-}\dif}^{+}(V)$, we obtain
\begin{align}\beta_{i}(t^{M}g)=t^{M}g+\frac{(c_{i}(\beta_{i}))^{1}}{1!}(\nabla^{(i)})^{1}(t^{M}g)+\frac{(c_{i}(\beta_{i}))^{2}}{2!}(\nabla^{(i)})^{2}(t^{M}g)+\cdots \end{align}
mod $(t,t_{1},\ldots,t_{e})^{r}D_{e\text{-}\dif}^{+}(V)$ for all $\beta_{i}\in\Gamma_{i}^{t^{M}g,r}$. Note that this series is a finite sum mod $(t,t_{1},\ldots,t_{e})^{r}D_{e\text{-}\dif}^{+}(V)$ by Proposition $3.5$. Thus, there exists  $L\in\mathbb{N}$ such that we have $(\nabla^{(i)})^{L}(t^{M}g)\not=0$ and $(\nabla^{(i)})^{L+1}(t^{M}g)=0$ mod $(t,t_{1},\ldots,t_{e})^{r}D_{e\text{-}\dif}^{+}(V)$.
On the other hand, since we have $(\nabla^{(i)})^{j}(g)\in (B_{\dR,K^{\pf}}\otimes_{\mathbb{Q}_{p}}V)^{H}$ by Lemma $3.3$ and $\nabla^{(0)}(\frac{1}{t^{j}}(\nabla^{(i)})^{j}(g))=0$ by Proposition $3.5$, there exists a finite field extension $M^{\pf}/K^{\pf}$ in $K_{\infty}^{\pf}$ such that $\verb+{+\frac{1}{t^{j}}(\nabla^{(i)})^{j}(g)\verb+}+_{j=0}^{L}$  is contained in $D_{\dR,M^{\pf}}(V)$. Write  $M^{\pf}=M_{0}^{\pf}(b)$. Then, since $V$ is also a crystalline representation of $G_{M^{\pf}}$,  
we have the equality $D_{\dR,M^{\pf}}(V)=M^{\pf}_{0}(b)\otimes _{M_{0}^{\pf}}D_{\cris, M^{\pf}}(V)$. Thus, we can write  
\begin{align}\frac{1}{t^{j}}(\nabla^{(i)})^{j}(g)=\sum_{m,n}b^{m}\otimes a_{ijmn}g_n \quad (a_{ijmn}\in \widetilde{B}_{\cris,K^{\pf}}).\end{align}
By $(3.3)$ and $(3.4)$, we obtain 
\begin{align}\beta_{i}(t^{M}g)=t^{M}\sum_{m,n}b^{m}\otimes(\sum_{j=0}^{L}\frac{(c_{i}(\beta_{i}))^{j}}{j!} a_{ijmn}t^{j})g_{n} \quad \bigl(\mod \ (t,t_{1},\ldots,t_{e})^{r}\bigr).\end{align}
On the other hand, since $\widetilde{D}_{\cris,K^{\pf}}(V)$ is stable under the action of $\Gamma_{i}^{t^{M}g,r}$ and  $\verb+{+b^{m}\verb+}+_{m=0}^{\delta-1}$ ($\delta=[M^{\pf}:M^{\pf}_{0}]$) is linearly independent over $\widetilde{B}_{\cris,K^{\pf}}$, the terms of the RHS of $(3.5)$ have to be $0$ for $m\not=0$. 
Then, for $m\not=0$, we have  $\sum_{j=0}^{L}\frac{\lambda^{j}}{j!} a_{ijmn}t^{j}=0$ for $\lambda$ in an open subgroup of $\mathbb{Z}_{p}$: this implies that $a_{ijmn}=0$ for $m\not=0$. In particular, we obtain $\nabla^{(i)}(g)\in \widetilde{D}_{\cris,K^{\pf}}(V)$ ($i\not=0$) by $(3.4)$. 
\end{proof}
\section{Proof of the main theorem}
In this section, we will give proofs only in the crystalline representation case since the semi-stable representation case is similar.
\begin{prop} 
We have the following implications.
\begin{enumerate}
\item If $V$ is a crystalline representation of $G_{K}$, then it is a crystalline representation of $G_{K^{\pf}}$.
\item If $V$ is a semi-stable representation of $G_{K}$, then it is a semi-stable representation of $G_{K^{\pf}}$.  
\end{enumerate}
\end{prop} 
\begin{proof}
Since $V$ is a crystalline representation of $G_{K}$, there exists a $G_{K}$-equivariant isomorphism of $B_{\cris,K}$-modules
\begin{align}B_{\cris,K}\otimes_{\mathbb{Q}_{p}}V\simeq (B_{\cris,K})^{d} \quad (d=\dim_{\mathbb{Q}_{p}}V).\end{align}
By tensoring (4.1) by $B_{\cris,K^{\pf}}$ over $B_{\cris,K}$ (which is induced by the $G_{K^{\pf}}$-equivariant surjection $p:B_{\cris,K}\twoheadrightarrow B_{\cris,K^{\pf}}:r_{i}\mapsto 0$),
we obtain a $G_{K^{\pf}}$-equivariant isomorphism of $B_{\cris,K^{\pf}}$-modules
$$B_{\cris,K^{\pf}}\otimes_{\mathbb{Q}_{p}}V\simeq (B_{\cris ,K^{\pf}})^{d}.$$
This means that $V$ is a crystalline representation of $G_{K^{\pf}}$.
\end{proof}
\begin{proof}[Proof of Theorem 1.1.] 
It remains to show that, if $V$ is a $p$-adic representation of $G_{K}$ whose  restriction to $G_{K^{\pf}}$ is crystalline, then $V$ is a potentially crystalline representation of $G_{K}$. 
Since $V$ is a crystalline representation of $G_{K^{\pf}}$, there exists a basis $\verb+{+g_{j}\verb+}+_{j=1}^{d}$ of $\widetilde{D}_{\cris,K^{\pf}}(V)$ over $\widetilde{B}_{\cris,K^{\pf}}$ which satisfies the properties in Proposition $3.4$. From this $\verb+{+ g_{j}\verb+}+_{j=1}^{d}$, for all finite extension $L/K$ in $\overline{K}$, we shall construct $L_{0}^{\pf}$-linearly independent elements $\verb+{+f_{j}\verb+}+_{j=1}^{d}$ in $B_{\cris,K}\otimes_{\mathbb{Q}_{p}}V$  
such that $\nabla^{(i)}(f_{j})=0$ for all $0\leq i\leq e$ and $1\leq j\leq d$. 

{\bf (A) Construction of $\verb+{+f_{j}\verb+}+_{j=1}^{d}$ in $B_{\cris,K}\otimes_{\mathbb{Q}_{p}}V$.} 

By Propositions $3.5$ and $3.6$, we have 
$(\nabla^{(1)})^{k_{1}}\cdots(\nabla^{(e)})^{k_{e}}(g_{j})=t^{k_{1}+\cdots +k_{e}}\sum_{l=1}^{d}$ $c_{j\underline{k}l}g_{l}$
where the $c_{j\underline{k}l}$ $(\underline{k}=(k_{1},\ldots,k_{e}))$ are elements of $B_{\cris,K^{\pf}}^{+}$  such that $\nabla^{(0)}(c_{j\underline{k}l})$ $=0$.
On the other hand, for $N\in\mathbb{N}$, we obtain,  
 \begin{align} \varphi^{N+1}((\nabla^{(1)})^{k_{1}}\cdots(\nabla^{(e)})^{k_{e}}(g_{j})) 
=(pt)^{k_{1}+\cdots+k_{e}}\sum_{l=1}^{d}p^{N(k_{1}+\cdots+k_{e})}\varphi^{N+1}(c_{j\underline{k}l}g_{l})\end{align}
where the $\varphi^{N+1}(c_{j\underline{k}l})$ are elements of $B_{\cris,K^{\pf}}^{+}$ such that  $\nabla^{(0)}(\varphi^{N+1}(c_{j\underline{k}l}))=0$.
Let $U_{i}$ denote the matrix which represents the action of $\nabla^{(i)}/t$ ($1\leq i\leq e$) with respect to the basis $\verb+{+g_{j}\verb+}+_{j=1}^{d}$ and take $N$ large enough such that we have $p^{N}U_{i}\in\M_{d}(A_{\cris,K^{\pf}})$ for all $1\leq i\leq e$. On the other hand, by applying the same method as in Proposition $3.6$ to the entries of $U_{i}$,  we have $\frac{\nabla^{(j)}}{t}(p^{N}U_{i})\in\M_{d}(\widetilde{B}_{\cris,K^{\pf}})$. Since we have $\nabla^{(0)}(\frac{\nabla^{(j)}}{t}(p^{N}U_{i}))=0$, this means that we obtain $\frac{\nabla^{(j)}}{t}(p^{N}U_{i})\in\M_{d}(L^{\pf}_{0})$ for a finite extension $L/K$ in $\overline{K}$ and, in particular,   $\frac{\nabla^{(j)}}{t}(p^{N}U_{i})\in\M_{d}(B_{\cris,K^{\pf}}^{+})$ ($1\leq i,j\leq e$). Furthermore, since $\nabla^{(j)}$ is the form $t\frac{d}{dt_{j}}$ on $K_{\infty}^{(\pf)}[[t,t_{1},\ldots,t_{e}]]$ and $\frac{\nabla^{(j)}}{t}=\frac{d}{dt_{j}}$ does not decrease the $p$-adic valuation of an element of  $K_{\infty}^{(\pf)}[[t,t_{1},\ldots,t_{e}]]\cap A_{\cris,K^{\pf}}\ (\subset B_{\dR,K})$, we obtain $\frac{\nabla^{(j)}}{t}(p^{N}U_{i})\in\M_{d}(A_{\cris,K^{\pf}})$ ($1\leq i,j\leq e$). Thus, it follows that we have $p^{N(k_{1}+\cdots+k_{e})}c_{j\underline{k}l}\in A_{\cris,K^{\pf}}$ and    $p^{N(k_{1}+\cdots+k_{e})}\varphi^{N+1}(c_{j\underline{k}l})\in  A_{\cris,K^{\pf}}$. Define $\verb+{+f_{j}\verb+}+_{j=1}^{d}\subset   B_{\cris,K}\otimes_{\mathbb{Q}_{p}}V$  by
$$f_{j}=\sum _{0\leq k_{1},\ldots,k_{e}}(-1)^{k_{1}+\cdots+k_{e}}\frac{t_{1}^{k_{1}}\cdots t_{e}^{k_{e}}}{k_{1}!\cdots k_{e}!t^{k_{1}+\cdots +k_{e}}}\varphi^{N+1}((\nabla^{(1)})^{k_{1}}\cdots (\nabla^{(e)})^{k_{e}}(g_{j}))$$
where $t_{i}=\log([\tilde{b_{i}}]/b_{i})$ denotes the element of $\Ker (\theta_{K_{0}}) \ (\subset A_{\cris,K})$.  
Note that this series converges in $B_{\cris,K}\otimes_{\mathbb{Q}_{p}}V$ for the $p$-adic topology by (4.2) and thus $f_{j}$ actually defines an element of $B_{\cris,K}\otimes_{\mathbb{Q}_{p}}V$. 
Then, it is easy to verify that we have $\nabla^{(i)}(f_{j})=0$ for all $1\leq i\leq e$ and $1\leq j\leq d$ by using the Leibniz rule. Furthermore, by using $(4.2)$ and the fact $\nabla^{(0)}(\varphi^{N+1}(g_{j}))=0$, we
can deduce that we have $\nabla^{(0)}(f_{j})=0$ for all $1 \leq  j \leq d.$  

{\bf (B) $\verb+{+f_{j}\verb+}+_{j=1}^{d}\subset B_{\cris,K}\otimes_{\mathbb{Q}_{p}}V$ is linearly independent over $L_{0}^{\pf}$.}

The homomorphism $p:B_{\cris,K}\twoheadrightarrow B_{\cris,K^{\pf}}$ induces
$$B_{\cris,K}\otimes_{\mathbb{Q}_{p}}V\twoheadrightarrow B_{\cris,K^{\pf}}\otimes_{\mathbb{Q}_{p}}V:f_{j}\mapsto \varphi^{N+1}(g_{j}).$$
Since $\verb+{+g_{j}\verb+}+_{j=1}^{d}$ forms a  basis of   $\widetilde{D}_{\cris,K^{\pf}}(V)$ over $\widetilde{B}_{\cris,K^{\pf}}$ and satisfies $\nabla^{(0)}(g_{j})=0$,  there exists a finite field extension $M/K$ in $\overline{K}$ such that $\verb+{+g_{j}\verb+}+_{j=1}^{d}$ forms a basis of $D_{\cris,M^{\pf}}(V)$ over $M_{0}^{\pf}$. Furthermore, since $\varphi: D_{\cris,M^{\pf}}(V)\rightarrow D_{\cris,M^{\pf}}(V)$ is bijective,   $\verb+{+\varphi^{N+1}(g_{j})\verb+}+_{j=1}^{d}$ also forms a basis of $D_{\cris,M^{\pf}}(V)$  over $M_{0}^{\pf}$. Thus, it follows  that $\verb+{+f_{j}\verb+}+_{j=1}^{d}$ is linearly independent over $L_{0}^{\pf}$ in $B_{\cris,K}\otimes_{\mathbb{Q}_{p}}V$ for all finite extension $L/K$ in  $\overline{K}$.


{\bf (C) Conclusion.}


The maps $\log(\gamma)$ and $\log(\beta_{i})$ ($1\leq i\leq e)$ act trivially on the $K_{0}$-vector space generated by $\verb+{+f_{j}\verb+}+_{j=1}^{d}$ (because $\nabla^{(0)}(f_{j})=0$ and $\nabla^{(i)}(f_{j})=0$ for all $1\leq i\leq e$ and $1\leq j\leq d$). This means that $\Gamma_{K}$ acts on this $K_{0}$-vector space via finite quotient and there exists a finite field extension $L/K$ in $\overline{K}$ such that $\verb+{+f_{j}\verb+}+_{j=1}^{d}$ forms a basis of $D_{\cris,L}(V)$ over $L_{0}\ (\subset L_{0}^{\pf}$).
     
\end{proof}
\section{The $p$-adic monodromy theorem of Fontaine in the imperfect residue field case}
In this section, we generalize the $p$-adic monodromy theorem of Fontaine to the imperfect residue field case.
Now, we recall the results of [Be] and [M].
\begin{thm} \verb+[+Be, Corollary 5.22.\verb+]+  
Let $L$ be a complete discrete valuation field of characteristic $0$ with perfect residue field of characteristic $p>0$ and $V$ be a $p$-adic representation of $G_{L}$. 
Then, $V$ is a de Rham representation of $G_{L}$ if and only if $V$ is a potentially semi-stable representation of $G_{L}$.
\end{thm}     
\begin{thm} \verb+[+M, Theorem 4.10.\verb+]+ 
Let $K$ be a complete discrete valuation field of characteristic $0$ with residue field $k$ of characteristic $p>0$ such that $[k:k^p]=p^{e}<\infty$ and $V$ be a $p$-adic representation of $G_{K}$. Let $K^{\pf}$ be the field extension of $K$ defined as before. Then,
$V$ is a de Rham representation of $G_{K}$ if and only if $V$ is a de Rham representation of $G_{K^{\pf}}$.
\end{thm}  
Since $K^{\pf}$ has  perfect residue field, we can apply Theorem $5.1$ to  the restriction of $V$ to  $G_{K^{\pf}}$.  
\begin{proof}[Proof of Corollary 1.2.] 
$V$ is a de Rham representation of $G_{K}$ if and only if $V$ is a de Rham representation of $G_{K^{\pf}}$ by Theorem $5.2$. Next, $V$ is a de Rham representation of $G_{K^{\pf}}$ if and only if $V$ is a potentially semi-stable representation of $G_{K^{\pf}}$ by Theorem $5.1$.
Finally, $V$ is a potentially semi-stable representation of $G_{K^{\pf}}$ if and only if $V$ is a potentially semi-stable representation of $G_{K}$ by  Theorem $1.1.$
\end{proof}

\end{document}